\documentclass[a4paper,12pt]{amsart}
\usepackage{amsmath,amssymb,latexsym,amsthm}

\newcommand\E{\mathcal E}
\newcommand\Ii{\mathcal I}

\newcommand\gH{\mathcal H}
\newcommand\Oc{\mathcal O}

\newcommand\Pj{\mathbb{P}}
\newcommand\ZZ{\mathbb{Z}}

\newtheorem{theorem}{Theorem}[section]
\newtheorem{proposition}[theorem]{Proposition}
\newtheorem{case}[theorem]{Case}
\newtheorem{definition}[theorem]{Definition}

\newtheorem{corollary}[theorem]{Corollary}

\newtheorem{remark}[theorem]{Remark}
\newtheorem{example}[theorem]{Example}
\newtheorem{lemma}[theorem]{Lemma}
\newtheorem{claim}[theorem]{Claim}
\theoremstyle{plain}
\theoremstyle{definition}
\theoremstyle{remark}
\begin{document}
\title{ A splitting criterion for rank $2$ bundles on a general sextic threefold }
\author{Luca Chiantini, Carlo Madonna}

\maketitle
\begin{abstract}
In this paper we show that on a general sextic hypersurface $X\subset\Pj^4$, a rank 2 vector bundle $\E$ splits if and only if $h^1(\E(n))=0$ for any $n \in \ZZ$. We get thus a characterization of complete intersection curves in $X$.\end{abstract}

\section{Introduction} \label{intro}

\indent Curves and vector bundles defined on a smooth projective threefold $X\subset\Pj^n$ have been considered as a main tool for the description of the geometry of $X$. Indeed, as soon as $X$ is sufficiently positive (e.g. Calabi-Yau or of general type) then one expects to have few types of curves and  bundles on it, so that these objects may work as sensible invariants for the threefold.\par

Clearly one obtains curves on $X$ just intersecting it with two hypersurfaces of $\Pj^n$, but it is a general non--sense that threefolds of general type should not contain many other curves. In fact, even in the most familiar case of general  hypersurfaces of high degree $r$ in $\Pj^4$, we do not know about the existence of curves in $X$ whose degree is not a multiple of $r$ (see e.g. \cite{Voi},\cite{Am}, \cite{Wu}, \cite{FR}). Similarly, one gets vector bundles of rank $2$ taking the direct sum of two line bundles on $X$, but it seems hard to find further examples on threefolds of general type.\par

The link between rank $2$ bundles and curves in smooth threefolds relies on the notion of subcanonical variety via Serre's correspondence (see e.g. \cite{Har}). Remind that a projective, locally Cohen--Macaulay variety $Y$ is {\it subcanonical}  when its dualizing sheaf $\omega_Y$ is $\Oc_Y(eH)$ for some integer $e$, $H$ being the class of a hyperplane section. If the threefold $X$ is subcanonical itself and moreover $h^1(\Oc_X(m))=h^2(\Oc_X(m))=0$ for all $m$ (notice that complete intersections satisfy  these conditions), then subcanonical curves $C$ on $X$ are exactly those curves which arise as zero--loci of global sections of a rank $2$ bundle $\E$ on $X$. There is a natural exact sequence:
$$
0\to \Oc_X \to \E\to \Ii_C(c_1(\E))\to 0 \eqno {(0)}
$$
where $\Ii_C$ is the ideal sheaf of the curve in $X$. Moreover $C$ is complete intersection {\it in $X$} if and only if the associated bundle $\E$ splits as a sum of line bundles.\par

The problem of determining conditions under which a curve is complete intersection was  studied by classical geometers. It is still open in higher dimensional projective spaces, but in $\Pj^3$ a solution was given by G. Gherardelli in 1942:

\begin{theorem} (Gherardelli, \cite {Ghe}) $C\subset \Pj^3$ is complete intersection
if and only if it is subcanonical and arithmetically Cohen-Macaulay.
\end{theorem}

In modern language, arithmetically Cohen--Macaulay means that the ideal sheaf $\Ii_C$ of $C$ satisfies $H^1(\Ii_C(n))=0$ for all $n$. By sequence (0), for subcanonical curves this condition is equivalent to the vanishing of all the cohomology groups $H^1(\E(n))$ of the associated rank $2$ bundle. We call rank $2$ bundles which satisfy this cohomological condition ''{\it arithmetically Cohen-Macaulay (ACM) bundles}'' (see Definition \ref{def:ACM}). Gherardelli's result was rephrased and generalized in the language of bundles by Horrocks:

\begin{theorem} (Horrocks, \cite {Ho}) A rank $2$ vector bundle on $\Pj^3$ splits in a sum of line bundles if and only if it is ACM.  
\end{theorem}

Gherardelli's theorem and Horrocks' criterion fail when $\Pj^3$ is replaced with a more general threefold. Even in the case of quadrics in $\Pj^4$, spinor vector bundles associated to lines are counterexamples (see \cite{Ott}). ACM bundles on some (mainly Fano) threefolds are studied in  the recent literature, and in some cases their moduli spaces are described. We refer to \cite{A-C} \cite{Bea2} \cite{M-T} and \cite{Dru} for  cubic threefolds, to \cite{I-M} and \cite{Mad2} for quartic threefolds and more  generally to \cite {A-G}, \cite{Mad3}, \cite{Kno}, \cite{Bea1}, \cite{S-W}, \cite{Kley}, \cite{Tho}  and \cite{B-G-S}. Some ACM bundles on Fano hypersurfaces of $\Pj^4$ are related to a pfaffian description of forms (see \cite {Bea1}). In a previous paper we proved that all stable ACM bundles, on a smooth quintic (Calabi-Yau) threefold, are rigid (\cite{CM}), giving a partial answer to a conjecture of Tyurin  (\cite{Tyu}) (see also \cite{Mad4}). \par

Some of the previous classification results are obtained by means of a theorem of the second author (\cite{Mad1}), who proved that Horrocks' splitting criterion works even for bundles $\E$ on smooth hypersurfaces in $\Pj^4$, under some numerical conditions on the invariants of $\E$ (see \ref{Mad} below).  It turns out that if the indecomposable ACM bundle $\E$ is normalized so that $h^0(\E)>h^0(\E(-1))=0$, then only few possibilities are left for its Chern classes.\par 

Using this reduction, we explore here the existence of ACM bundles and subcanonical curves on hypersurfaces $X$ of general type in $\Pj^4$. Of course, when $X$ is an arbitrary hypersurface, we cannot expect to say much about its ACM bundles. Conversely, when $X$ is general, one may hope to control the situation. \par

The main result of this paper concerns general sextic threefolds. We prove that indeed  a general sextic $X$ has no indecomposable rank $2$ ACM bundles; in other words Horrocks criterion works for rank $2$ bundles on a general sextic:   

\begin{theorem} \label{thm:splitting} Let $\E$ be a rank 2 vector bundle
on a general sextic threefold $X$. Then $\E$ splits if and only 
if $H^1(\E(n))=0$ for all $n \in \ZZ$.
\end{theorem}

Using Serre's correspondence, the result can be rephrased for curves in $X$ in the following, Gherardelli's type criterion:

\begin{corollary} A (locally complete intersection) curve $C$ contained in a general sextic hypersurface $X\subset\Pj^4$ is complete intersection in $X$ if and only if it is subcanonical and arithmetically Cohen--Macaulay.
\end{corollary}

Notice that none of the two assumptions of the previous corollary can be dropped:

\begin{example} There are ACM curves on a general sextic threefold, which are not subcanonical (hence are not complete intersection).\par\noindent
\rm An example was found by Voisin in \cite{Voi}, starting with $2$ plane sections of $X$ which meet at a point, and using linkage. In these examples the degree is always a multiple of 6.
\end{example}

\begin{example} It is not hard to find examples of smooth irreducible subcanonical curves in a general  sextic threefold, which are not ACM (hence not complete intersection).\par\noindent
\rm Just take two disjoint plane sections of $X$. Their union $Y$ is a subcanonical curve which is not ACM. If $\E$ is a rank $2$ bundle associated to $Y$, then the zero--locus of a general section of $\E(k)$, $k\gg 0$, has the required properties.\end{example}

The proof of our main result is achieved first using \ref{Mad} to get a rough classification of curves arising as zero--loci of ACM bundles on a sextic threefold (see section 3). 
Then, as in \cite{CM}, we use the method introduced by Kleppe and Mir\'o--Roig in \cite{K-M} to understand the infinitesimal deformations of the corresponding ACM subcanonical curves, and we get rid of all the possibilities. The splitting criterion is then established by a case by case analysis.   \par

As the number of cases which cannot be ruled out directly by \ref{Mad} increases with the degree of $X$, an extension of this procedure to hypersurfaces of higher degree looks unreliable.
We wonder if, based on our result, a degeneration argument could prove Horrocks' splitting criterion for general hypersurfaces of degree bigger than $6$. Also we would know the geometry of the variety of sextics which have some indecomposable ACM bundles with given Chern classes (somehow an analogue of Noether--Lefschetz loci for surfaces in $\Pj^3$).\par 

We finally observe that several numerical refinements of the main result, in the spirit of \cite{CV}, are immediate using our main theorem and theorem 3.8 of \cite{Mad1}. For instance one gets:

\begin{corollary}
Let $\E$ be a rank 2 vector bundle on a general sextic threefold
$X$. Then $\E$ splits if and only if
$h^1(\E(a))=0$, where $a=\frac{-c_1+3}2$ if $c_1$ is odd and
$a=\frac{-c_1+2}2$ if $c_1$ is even.
\end{corollary}

\section{Generalities} \label{S:gener} \par
We work in the projective space $\Pj^4$ over the complex field. 
We will denote by $\Oc$ the structure sheaf of $\Pj^4$.\par\smallskip
Let $X$ be a  general hypersurface of degree $r\geq 3$ in $\Pj^4$.
$X$ is smooth and we identify its Picard group with $\ZZ$, generated by the class
of an hyperplane section.
We use this isomorphism to identify  line bundles with integers. 
In particular, for any vector bundle $\E$ on $X$, we consider 
$c_1(\E) \in \ZZ$ and we write $\E(n)$ for 
$\E \otimes \Oc_X(n)$. 
When $\E$ has rank $2$, we have the following formulas for the Chern classes of the twistings of $\E$:
\[
\begin{aligned}
c_1(\E(n)) & = c_1(\E)+2n \\
c_2(\E(n)) & = c_2(\E)+rnc_1(\E)+rn^2.
\end{aligned}
\]
 
\begin{definition} \label{def:ACM}
Let $\E$ be a rank $2$ vector bundle on
a smooth projective threefold $X\subset \Pj^4$. We say that 
$\E$ is an {\it arithmetically Cohen--Macaulay (ACM for short) bundle} 
if $h^i(\E(n))=0$ for  $i=1,2$ and for any $n \in \ZZ$.
\end{definition}

Let us define the number:
\[
b(\E)=b:=\max \{ n \mid h^0(\E(-n)) \ne 0 \}.
\]

\begin{definition} 
We say that the rank $2$ bundle $\E$ is {\it normalized} when $b=0$. 
\end{definition}

Of course, after replacing $\E$ with the twist $\E(-b)$, 
we may always assume that it is normalized. 
Since the Picard group of $X$ is generated over $\ZZ$ by the 
hyperplane class, we get: 

\begin{remark} \label{stab} \rm
A rank 2 vector bundle $\E$ is semi--stable if and only if 
$2b-c_1\leq 0$. It is stable if and only if the strict inequality holds.
\par\noindent
The number $2b-c_1$ is invariant by twisting i.e. for all $n\in\ZZ$:
\[
2b-c_1=2b(\E(n))-c_1(\E(n)).
\]
This invariant measures the {\it level of stability} of $\E$.
\end{remark}

\noindent If $b=0$, it is shown in \cite{Har} 
Remark 1.0.1 that $\E$ has some 
global section whose zero--locus $C$ has codimension $2$. 
$C$ is a {\it subcanonical} curve of degree $c_2(\E)$, 
whose canonical divisor is $\omega_C=\Oc_C(c_1(\E)+r-5)$.\smallskip

\begin{remark} \label {rr} \rm Let $\E$ be a rank 2 vector bundle on $X$.
Since $\omega_X=\Oc_X(r-5)$, Serre's duality says that:
$$ 
h^3(\E(n)) = h^0(\E^\vee(-n+r-5))= h^0(\E(-c_1-n+r-5)).
$$
Moreover one computes:
$$
\begin{aligned}
\chi(\E) & =\frac{rc_1^3}6+\frac{5-r}4rc_1^2-\frac{5-r}2c_2 
-\frac{c_1c_2}{2} + \\
& +\frac{rc_1}{12}(2r^2-15r+35)+\frac r{12}(-r^3+10r^2-35r+50). 
\end{aligned}
$$
\end{remark}

\begin{remark} \label {6rr} \rm
When $X$ has degree $6$, the previous formula reduces to:
$$
\chi(\E)=c_1^3-\frac32c_1^2+\frac{c_2}2-\frac{c_1c_2}{2}+\frac{17}2c_1-8. 
$$
\end{remark}

\noindent We are going to use the main result of \cite{Mad1}, to get rid 
of most values of $c_1$ for non-splitting rank $2$ 
ACM bundles on smooth hypersurfaces of $\Pj^4$. We recall the result  here:

\begin{theorem} \label{Mad}
Let $\E$ be a normalized rank 2 ACM bundle on a 
smooth hypersurface $X\subset\Pj^4$ of degree $r$. Then, if $\E$ is indecomposable: 
$$2-r<c_1(\E)<r.$$ 
\end{theorem}

\section{ACM bundles with small $c_1$}

When the first Chern class of the {\it normalized} rank $2$ bundle $\E$ on $X$ is smaller than or equal to $6-r$, then $\E$ has a  section whose zero--locus is a curve $C$ with canonical class $\omega_C=\Oc_C(e)$, and $e=c_1(\E)+r-5\leq 1$.\par
Of course, we do not know much about $C$: it can be reducible, non reduced. We just know that it is locally complete intersection and subcanonical (observe that, in particular, its canonical sheaf is locally free). On the other hand, when $\E$ is ACM, for these cases of  low $c_1$ we can give a description of the invariants of $C$ 
(and sometimes of $C$ itself). Let us study this description case by case.
\smallskip

{\it Let $\E$ be a normalized rank 2 ACM 
bundle on a smooth threefold $X\subset\Pj^4$ of degree $r$.}

\begin{case}\label{retta}
Assume that $c_1(\E)=3-r$. Then $c_2(\E)=1$ and $\E$ has a 
section whose zero--locus is a line.
\end{case}
\begin{proof} $\E$ has a section whose zero--locus $C$ is a curve. 
The exact sequence (0) of the introduction here reads:
$$0\to \Oc_X \to \E\to \Ii_C(3-r)\to 0.$$
Consider $\E'=\E(r-3)$. $\E'$ is ACM, so that in particular $h^1(\E')=h^2(\E')=0$. Furthermore by the previous sequence $h^0(\E')=h^0(\Oc_X(r-3))$ and $h^3(\E')=h^3(\E(r-3))= h^0(\E(r-5))=h^0(\Oc_X(r-5))$. 
So one may compute $\chi(\E')$ directly.
On the other hand $\E'$ has Chern classes $c_1(\E')=r-3$ and $c_2(\E')=c_2(\E)$, and one can compute $\chi(\E')$ using the Riemann-Roch formula of remark \ref{rr}. 
After some easy computations, it turns out that $r$ disappears and one gets $c_2(\E)=1$, i.e. $C$ has degree 1.
\end{proof}

\begin{case}\label{conic}
Assume that $c_1(\E)=4-r$. Then $c_2(\E)=2$ and $\E$ has a section whose zero--locus is a (possibly singular) conic.
\end{case}
\begin{proof} As above, $\E$ has a section whose zero--locus is a curve; call it $C$ 
and consider the exact sequence:
$$0\to \Oc_X \to \E\to \Ii_C(4-r)\to 0.$$
Take $\E'=\E(r-3)$. $\E'$ is ACM, so that in particular $h^1(\E')=h^2(\E')=0$. 
Furthermore the previous sequence says that $h^0(\E')=h^0(\Oc_X(r-3))+h^0(\Ii_C(1))$ 
and 
$h^3(\E')= h^0(\E(r-6))=h^0(\Oc_X(r-6))$. 
So one can compute $\chi(\E')$ in terms of $r$ and $h^0(\Ii_C(1))$.
On the other hand $\E'$ has Chern classes $c_1(\E')=r-2$ and $c_2(\E')=c_2(\E)+r^2-3r$, 
and one can compute $\chi(\E')$ using the Riemann-Roch formula of remark \ref{rr}.
After some easy computations, it turns out that $r$ disappears and one gets:
$$\frac 32 c_2(\E)=5-h^0(\Ii_C(1)).$$
Since $c_2(\E)=\deg(C)$ is a positive integer, the unique possibility is 
$c_2(\E)=2$ and $h^0(\Ii_C(1))=2$, i.e. $C$ is a plane curve of degree $2$.
\end{proof}

\begin{case}\label{elliptic}
Assume that $c_1(\E)=5-r$. Then only the following possibilities 
may arise:
\begin{enumerate}
\item $c_2(\E)=3$ and $\E$ has a section whose zero--locus is a plane cubic;
\item $c_2(\E)=4$ and $\E$ has a section whose zero--locus is a complete intersection elliptic space curve;
\item $c_2(\E)=5$ and $\E$ has a section whose zero--locus is a non--degenerate elliptic curve.
\end{enumerate}
\end{case}
\begin{proof} As above we know that $\E$ has a section whose zero--locus $C$ is a curve, with an exact sequence:
$$0\to \Oc_X \to \E\to \Ii_C(5-r)\to 0.$$
Take here $\E'=\E(r-4)$. As usual $h^1(\E')=h^2(\E')=0$, while the previous sequence says that 
$h^0(\E')=h^0(\Oc_X(r-4))+h^0(\Ii_C(1))$ and $h^3(\E')= h^0(\E(r-6))=h^0(\Oc_X(r-6))$. So, as above, 
one can compute $\chi(\E')$ in terms of $r$ and $h^0(\Ii_C(1))$.
On the other hand $\E'$ has Chern classes $c_1(\E')=r-3$ and $c_2(\E')=c_2(\E)+r^2-4r$, and one 
can compute $\chi(\E')$ 
using the Riemann-Roch formula. 
It turns out that $r$ disappears and one gets:
$$c_2(\E)=5-h^0(\Ii_C(1)).$$
Observe that $c_2(\E)<2$ is impossible, since otherwise the curve $C$ would be contained in too many independent hyperplanes.\par\noindent
If $c_2(\E)=3$, then $h^0(\Ii_C(1))=2$ and $C$ is a plane cubic.\par\noindent
If $c_2(\E)=4$, then $h^0(\Ii_C(1))=1$ and $C$ is an arithmetically Cohen-Macaulay space curve. 
It is well known that any arithmetically Cohen--Macaulay subcanonical curve in $\Pj^3$ is complete 
intersection (see e.g. \cite{Har}). The invariants tell us then that $C$ is complete intersection of 
two quadrics in $\Pj^3$.\par\noindent
Finally when $c_2(\E)=5$, then $C$ is a non--degenerate, elliptic quintic in $\Pj^4$.
\end{proof}

\begin{case}\label{canonica}
Assume that $c_1(\E)=6-r$. Then only the following possibilities may
arise:
\begin{enumerate}
\item $c_2(\E)=4$ and $\E$ has a section whose zero--locus is a plane quartic;
\item $c_2(\E)=6$ and $\E$ has a section whose zero--locus is a complete intersection space curve, 
of type (2,3);
\item $c_2(\E)=8$ and $\E$ has a section whose zero--locus is a non--degenerate (possibly singular) 
canonical curve of genus 5.
\end{enumerate}
\end{case}
\begin{proof} Call $C$ the zero--locus of a non--trivial section of $\E$. 
The exact sequence (0) reads:
$$0\to \Oc_X \to \E\to \Ii_C (6-r)\to 0.$$
Take $\E'=\E(r-6)$. Since $h^1(\E')=h^2(\E')=0$, 
$h^0(\E')=h^0(\Oc_X(r-6))$ and $h^3(\E')=
h^0(\Oc_X(r-5))+h^0(\Ii_C(1))$, one can compute $\chi(\E')$ in 
terms of $r$ and $h^0(\Ii_C(1))$. $\E'$ has Chern classes $c_1(\E')=r-6$ and $c_2(\E')=c_2(\E)$, 
and one can also 
compute $\chi(\E')$ using the Riemann-Roch formula. One gets:
$$c_2(\E)=8-2h^0(\Ii_C(1)).$$
If $h^0(\Ii_C(1))>0$, then $C$ is degenerate, so it must be complete intersection, and 
one concludes exactly as above.\par\noindent
Otherwise $c_2(\E)=8$ and $\omega_C=\Oc_C(1)$ 
implies that $C$ has arithmetic genus 5. \par\noindent
\end{proof}

Observe that by Riemann-Roch and Serre's duality, in the previous example, when $\deg(C)=8$ then $h^0(\Ii_C(2))=3$. If the $3$ quadrics intersect along a curve $C$, 
then $C$ is complete intersection in $\Pj^4$.
\smallskip

Going further, we can only decribe $c_2(\E)$ in terms of the number 
of independent hypersurfaces of degree $2,3,\dots$ through $C$, 
and the number of possibilities grows considerably.\par
On the other hand, in the case $r=6$, which is the object 
of our study, sometimes it is  still possible to determine $c_2(\E)$ 
accurately.\par\noindent
The worst cases, for $r=6$, arise from $c_1(\E)=7-r$ and $c_1(\E)=8-r$.

\begin{case}\label{bicanonica}
Assume that $c_1(\E)=7-r$. Then $c_2(\E)\leq 14$. More precisely, $c_2(\E)=
14 - h^0(\Ii_C(2))$, $C$ being the zero--locus of a non--trivial section of $\E$. 
\end{case}
\begin{proof} Here sequence (0) reads:
$$0\to \Oc_X \to \E\to \Ii_C (7-r)\to 0.$$
Take $\E'=\E(r-5)$. Since $h^1(\E')=h^2(\E')=0$, $h^0(\E')=h^0(\Oc_X(r-5))+
h^0(\Ii_C(2))$ and $h^3(\E')=h^0(\Oc_X(r-7))$, one can compute $\chi(\E')$ 
in terms of $r$ and $h^0(\Ii_C(2))$. $\E'$ has Chern classes $c_1(\E')=r-3$ and $c_2(\E')=c_2(\E)+2r^2-10r$, 
and one can also compute $\chi(\E')$ using Riemann-Roch. 
Comparing the two computations, as usual, $r$ disappears and one gets:
$$c_2(\E)=14-h^0(\Ii_C(2)).$$
\end{proof}

\begin{case}\label{tricanonica}
Assume that $c_1(\E)=8-r$. Then $c_2(\E)\leq 20$. 
More precisely, $c_2(\E)=20 - 2h^0(\Ii_C(2))+2h^0(\Ii_C(1))$, $C$ 
being the zero--locus of a non--trivial section of $\E$. 
\end{case}
\begin{proof} Here sequence (0) reads:
$$0\to \Oc_X \to \E\to \Ii_C (8-r)\to 0.$$
Take $\E'=\E(r-6)$. Since $h^1(\E')=h^2(\E')=0$, 
$h^0(\E')=h^0(\Oc_X(r-6))+h^0(\Ii_C(2))$ and 
$h^3(\E')=h^0(\Oc_X(r-7))+h^0(\Ii_C(1))$, one can compute $\chi(\E')$ in terms 
of $r$, $h^0(\Ii_C(1))$ and $h^0(\Ii_C(2))$. One can also compute $\chi(\E')$ using Riemann-Roch. Comparing the two computations, one gets:
$$c_2(\E)=20-2h^0(\Ii_C(2))+2h^0(\Ii_C(1)).$$
\end{proof}

\begin{remark} \rm
In the two last cases, the degree of the curve $C$ (i.e. the second Chern class of $\E$) is in fact well--determined, unless $C$ belongs to some hyperplane or to some quadric hypersurface $Q$.\par
We know very well ACM subcanonical curves lying in hyperplanes: they are complete intersection by Horrocks' splitting criterion. \par
Also we know very well ACM subcanonical curves lying in \em{smooth} \rm quadrics. Observe however that $C$ itself might be reducible or non--reduced, so we have few informations on the singularities of $Q$. In particular, at this stage, we cannot use 
$Q$ to determine a classification for $C$. 
\end{remark}

The proof for the following cases works exactly as above and we omit it.

\begin{case}\label{quadricanonica}
Assume that $c_1(\E)=9-r$. Then $c_2(\E)=30 - h^0(\Ii_C(3))+h^0(\Ii_C(1))$, 
$C$ being the zero--locus of a non--trivial section of $\E$. 
\end{case}

\begin{case}\label{pentacanonica}
Assume that $c_1(\E)=10-r$. Then $c_2(\E)=40 - 2h^0(\Ii_C(3))+2h^0(\Ii_C(2))$, 
$C$ being the zero--locus of a non--trivial section of $\E$. 
\end{case}

\begin{case}\label{esacanonica}
Assume that $c_1(\E)=11-r$. Then $c_2(\E)=55 - h^0(\Ii_C(4))+h^0(\Ii_C(2))$, 
$C$ being the zero--locus of a non--trivial section of $\E$. 
\end{case}

The reader can easily find a generalization of the previous 
computations for every $c_1$.\par
We rephrase these last cases in our situation, where $r=6$.

\begin{proposition} With $X$, $\E$ as above, assume furthermore that $r=6$. Then:
\begin{enumerate}
\item if $c_1(\E)=2$ then $c_2(\E)=20 - 2h^0(\Ii_C(2))$;
\item if $c_1(\E)=3$ then $c_2(\E)=30 - h^0(\Ii_C(3))$;
\item if $c_1(\E)=4$ then $c_2(\E)=40$;
\item if $c_1(\E)=5$ then $c_2(\E)=55$.
\end{enumerate} 
\end{proposition}
\begin{proof} Consider $c_1(\E)=5=11-r$. Sequence (0) reads as:
$$0\to \Oc_X \to \E\to \Ii_C (5)\to 0$$
and since $\E$ is normalized, it turns out that $h^0(\Ii_C(4))=h^0(\Ii_C(2))=0$. 
The claim follows now from case \ref{esacanonica}.
A similar argument concludes the remaining cases.
\end{proof}

\begin{proposition} Assume $r=6$. Let $C$ be the zero--locus of a 
non--trivial section of $\E$. 
\begin{enumerate}
\item If $c_1(\E)=4$ (hence $c_2(\E)=40$) then the ideal of $C$ is generated in degree $\leq 5$;
\item if $c_1(\E)=5$ (hence $c_2(\E)=55$) then $C$ is smooth, irreducible and 
its ideal is generated by quintics.
\end{enumerate}
\end{proposition}
\begin{proof} Consider $c_1(\E)=5=11-r$. Sequence (0) reads as:
$$0\to \Oc_X \to \E\to \Ii_C (5)\to 0$$
One knows that $h^1(\E(-1))=h^2(\E(-2))=0$ while $h^3(\E(-3))=h^0(\E(-1))=0$. 
It follows that $\E$ is regular in the sense of Castelnuovo--Mumford, hence it 
is generated by global sections. It is well--known that, in this case, working in characteristic $0$, the zero--locus of a general section of $\E$ is smooth. It is also connected, since $C$ is ACM. Finally $\Ii_C(5)$, which is a quotient of $\E$, is also generated by global sections since $h^1(\Oc_X)=0$.\par\noindent
The case $c_1(\E)=4$ is similar: one shows that $\E(1)$ is regular.
\end{proof}

We summarize the previous discussion, for $r=6$, in the following table:

$$\begin{matrix}
c_1(\E) & c_2(\E) & informations \cr 
- 3 & 1 & \text{line} \cr
-2 & 2 & \text{conic} \cr
 -1 & \bigg\{ \begin{matrix} 3 \cr 4 \cr 5 \end{matrix} & \begin{matrix}
 \text{plane cubic} \cr
 \text{space curve c.i. type (2,2)} \cr
\text{elliptic non--degenerate} \end{matrix} \cr
0 & \bigg\{ \begin{matrix} 4 \cr 6 \cr 8 \end{matrix} & \begin{matrix}
 \text{plane quartic} \cr
 \text{space curve c.i. type (2,3)} \cr
\text{canonical non--degenerate} \end{matrix} \cr
1 & \leq 14 & h^0\Ii_C(2)=14-c_2(\E) \cr
2 & \leq 20 & \text{non--degenerate, }2h^0\Ii_C(2)=20-c_2(\E)\cr
3 & \leq 30 & h^0\Ii_C(3)=30-c_2(\E) \cr
4 & 40 & \text{generated by quintics} \cr
5 & 55 & \text{smooth, irreducible, generated by quintics} \cr
\end{matrix}$$ 
\bigskip

\section{Curves and sextic threefolds}

In order to prove the main result of the paper, it is sufficient to show that 
a general sextic threefold does not contain any curve of the types listed in 
the previous table. We do this with an examination, case by case, of the 
corresponding Hilbert scheme.
Let us denote by: \par\noindent
-- $\gH'_{c_1,c_2}$ the (locally closed) Hilbert scheme of 
ACM curves in $\Pj^4$ with degree $c_2$ and genus $1+(c_1+1)c_2/2$;\par\noindent
-- $\gH:=\gH_{c_1,c_2}$ the union of all the components of $\gH'_{c_1,c_2}$ 
containing points which correspond to curves, zero--loci of sections of non--splitting normalized rank 2 ACM bundles $\E$, with Chern classes $c_1(\E)$, $c_2(\E)$, defined over some smooth sextic threefold; \par\noindent
-- $N(C)$ the normal bundle in $\Pj^4$ of a curve $C$ in $\gH$. 
Of course at $C$:
$$\dim(\gH)\leq h^0(N(C));$$
-- $\Ii(C)$ (or $\Ii$ if no confusion arises) the ideal sheaf in $\Pj^4$ of a curve $C$.\par
Let $\Pj = \Pj^{209}$ be the scheme which 
parametrizes sextic threefolds in $\Pj^4$.
In the product $\gH\times\Pj$ one has the incidence 
variety (i.e. the {\it Hilbert flag scheme}) 
$$I = \{(C,X): X\text{ is smooth and } C\subset X\},$$
with the two obvious projections $p:I\to \gH$ and $q:I\to \Pj$.

\begin{proposition} \label{main}
Fix $c_1,c_2$ such that $\gH=\gH_{c_1,c_2}$ is not empty (in particular, $c_1\leq 5$). Assume that for any curve $C$, which is the zero--locus of a section of a normalized rank 2 ACM bundle on some smooth sextic, one has:
$$  h^0(\Ii(6))+h^0(N(C))-1<209 \eqno(1)$$
(recall that $\Ii$ denotes the ideal sheaf of $C$ in $\Pj^4$).\par\noindent
Then the map $q:I\to \Pj$ of the corresponding incidence variety $I$ to the 
parametrizing space of sextics, is not dominant.
\end{proposition}
\begin{proof} 
By construction any curve $Y$ in $\gH$ is ACM. In particular $h^1(\Ii(Y)(6))=0$. From the usual exact sequence one gets $210=h^0(\Oc(6))=h^0(\Ii(Y)(6))+h^0(\Oc_Y(6))$. Thus, by semicontinuity, both $h^0(\Ii(Y)(6))$ and $h^0(\Oc_Y(6))$ are constant in any fixed component of $\gH$, for their sum is constant. Hence $h^0(\Ii(Y)(6))$ does not depend on $Y$. Now observe that any component of $\gH$ contains an ACM $e$-subcanonical curve $C$, so we may use this curve to compute the number of independent sextic hypersurface passing through any $Y$ in $\gH$. Since $h^0(\Ii(Y)(6))-1$ is exactly the dimension of the fiber of $p:I\to\gH$ at $Y$, the claim follows by an easy dimensional count: indeed one gets $\dim(p^{-1}\gH)<209$.
\end{proof}

Our strategy now consists in bounding the dimension of $\gH$ 
by computing $h^0(N(C))$ for 
ACM curves of the types listed in the previous section.\smallskip

Curves arising as zero--loci of sections of rank 2 ACM bundles on $X$ 
are subcanonical ACM curves in $\Pj^4$. 
In particular, such curves are {\it arithmetically Gorenstein} (\cite{B-E}). The resolution of the 
ideal sheaf of arithmetically Gorenstein curves in $\Pj^4$ is described by the following:

\begin{proposition}\label{Gor}
Let $C$ be an $e$--subcanonical, ACM curve in $\Pj^4$ and call $\Ii$ the ideal sheaf of $C$ in $\Pj^4$. 
Then one has a resolution:
$$0\to \Oc(-e-5)\to \oplus_i \Oc(-b_i)  \to \oplus_i \Oc(-a_i)\to \Ii\to 0 \eqno{(2)}$$ 
which is self--dual, up to twisting. Hence if one orders the $a_i$'s and the $b_i$'s so 
that $a_1\leq\dots\leq a_n$ and $b_n\leq\dots\leq b_1$ (observe that the two orders are reversed), then:
$$ \forall i\quad -a_i = b_i-e-5.$$
\end{proposition}
\begin{proof} see \cite{B-E}, p.466.
\end{proof} 

Since a curve $C$ which is zero--locus of a section of a rank 2 ACM bundle on $X$ is locally complete intersection, then the (embedded) normal bundle of $C$ in $\Pj^4$ is well defined. The cohomology of this bundle is computed from the following formula  of Kleppe and Mir\'o--Roig (see \cite{K-M}): 

\begin{theorem}\label{KM}
Let $C$ be an $e$--subcanonical, ACM curve in $\Pj^4$. 
Then, with the previous notation, one can compute $h^0(N(C))$ from the formula:
$$ h^0(N(C))= 
\sum_{i=1}^n h^0(\Oc_C(a_i))+\sum_{1\leq  i<j\leq n} 
\binom{-a_i+b_j+4}4 +$$ 
$$-\sum_{1\leq i< j\leq n} \binom{a_i-b_j+4}4 - \sum_{i=1}^n 
\binom{a_i+4}4 \eqno{(3)}$$
\end{theorem}
\begin{proof} see \cite{K-M} th.2.6.
\end{proof} 

\section{ACM bundles on general sextic threefolds split}

\begin{case} There are no indecomposable rank $2$ ACM bundles $\E$ on a general 
sextic, with $c_1(\E)=-3$.
\end{case}
\begin{proof}
We know from 3.1 that $\E$ has sections which vanish on a line in $X$. So it is enough to remind that general sextic hypersurfaces contain no lines. Indeed the space of lines in $\Pj^4$ has dimension $6$ and any line is contained in a $202$-dimensional space of sextics, so the incidence 
variety $I$ has dimension $208$ here and it cannot dominate $\Pj$.  
\end{proof}

\begin{case} There are no indecomposable rank $2$ ACM bundles $\E$ on a general sextic, with $c_1(\E)=-2$.
\end{case}
\begin{proof}
Argue as above: by 3.2 the existence of $\E$ implies the existence of a conic in $X$. So it is enough to remind that general sextic hypersurfaces contain no conics. 
\end{proof}

\begin{case} There are no indecomposable rank 2 ACM bundles $\E$ on a general sextic, with $c_1(\E)=-1$.
\end{case}
\begin{proof}
We need to show that a general sextic hypersurface contains no 
elliptic ACM curves $C$ of degree $\leq 5$. This is well--known when 
$C$ is irreducible and reduced, but we need to show 
the non--existence in general. \par\noindent
Assume $\deg(C)=3$. Then $C$ is a plane cubic (hence complete 
intersection in $\Pj^4$), so its normal bundle in $\Pj^4$ has $h^0(N(C))=15$. 
Since $C$ is ACM, one computes $h^0(\Ii(6))=192$. Then just apply proposition 
\ref{main}.\par\noindent
If $\deg(C)=4$, then $C$ is still degenerate and complete intersection. One 
computes $h^0(N(C))=20$, $h^0(\Ii(6))=186$ and the claim follows.\par\noindent
Finally for $\deg(C)=5$ the curve is non--special, i.e. $h^1\Oc_C(1)=0$ for $C$ is ACM. Using  the Euler sequence, a computation yields $h^0(N(C))=25$ 
(see e.g. \cite{Kley}) and $h^0(\Ii(6))=180$: 
once again $I$ cannot dominate $\Pj$.
\end{proof}

\begin{case} There are no indecomposable rank 2 ACM bundles $\E$ on a general sextic, with $c_1(\E)=0$.
\end{case}
\begin{proof}
We want to show that a general sextic hypersurface contains no canonical ACM curves $C$ of (even) degree $\leq 8$. This is easy for $\deg(C)=4,6$, for the corresponding curves are degenerate, hence complete intersection, so one can easily argue as in the previous case. \par\noindent 
For $\deg(C)=8$, the situation is more delicate. If we know that the Hilbert scheme is irreducible, then we can still use complete intersection curves to determine $h^0(N(C))$. Unfortunately, we have few informations on $C$ and we cannot exclude, ''a priori'', the existence of a component of $\gH$ which contains (possibly reducible or non--reduced) canonical ACM curves, and whose general point is not complete intersection.\par\noindent
So instead, we compute $h^0(N(C))$ using theorem \ref{KM}, which indeed works also for ''bad'' curves.\par\noindent
Since $C$ is canonical and non--degenerate in our case, the three-steps free resolution
(2) reads:
$$0\to \Oc(-6)\to \Oc(-3)^x\oplus\Oc(-4)^3  \to \Oc(-3)^x\oplus\Oc(-2)^3\to \Ii\to 0 $$ 
Observe that indeed we cannot have quartics among the minimal generators otherwise, by auto--duality, we get a quadric syzygy, absurd since the resolution is minimal.\par\noindent
We do not know in principle how many independent cubics are among the minimal generators of $\Ii$. Nevertheless, when we apply formula (3), the contributions of $x$ cancel and $h^0(N(C))$ does not depend on it. One gets $h^0(N(C))=36$; since 
$h^0(\Ii(6))=166$, formula (4) shows that a general sextic does not contain curves like $C$. 
 \end{proof}

It will be true indeed in several cases that, even if we do not know some Betti 
number for the resolution of $\Ii$, nevertheless the contributions of the unknown generators cancel. \smallskip

In the previous cases, the rank $2$ bundles we worked with, were unstable or semi--stable. Let us turn our attention to the stable cases. We dispose first of the easiest situations, in which $c_1$ is big.

\begin{case} There are no indecomposable rank 2 ACM bundles $\E$ on a general sextic, with $c_1(\E)=5$.
\end{case}
\begin{proof}
We want to show that a general sextic hypersurface contains no $6$--subcanonical ACM curves $C$ of degree $55$ (and genus $166$).\par\noindent
We know that the ideal sheaf is generated by quintics. Furthermore by sequence (0) and since $\E$ is normalized, we know that $C$ lies in no quartics. One computes $h^0(\Oc_C(5))= h^0(\Oc_C(1))+110=115$. So by auto--duality, the free 
resolution (2) here reads:
$$0\to \Oc(-11)\to \Oc(-6)^{11}  \to \Oc(-5)^{11}\to \Ii\to 0 $$ 
Everything is determined and one can use (3) to bound the dimension of the Hilbert scheme. It turns out that $h^0(N(C))=154$ while $h^0(\Ii(6))=44$. Hence $\dim(I)\leq 197$ and by proposition 4.1 a general sextic does not contain curves like $C$. 
 \end{proof}

\begin{case} There are no indecomposable rank 2 ACM bundles $\E$ on a general sextic, with $c_1(\E)=4$.
\end{case}
\begin{proof}
We must exclude the existence of $5$--subcanonical ACM curves $C$ of degree $40$ (and genus $101$) on a general sextic hypersurface $X$.\par\noindent
We know from section 3 that the ideal sheaf of $C$ in $\Pj^4$ is generated by quartics and quintics, hence  by auto--duality, the free resolution (2) here reads:
$$0\to \Oc(-10)\to \Oc(-6)^5\oplus\Oc(-5)^x  \to \Oc(-5)^x\oplus\Oc(-4)^5\to \Ii\to 0 $$ 
Using (3) to compute $h^0(N(C))$, it turns out that the unknown $x$ cancels. 
One computes $h^0(N(C))=125$ while $h^0(\Ii(6))=70$. Hence here $\dim(I)\leq 194$ and a general sextic does not contain curves like $C$. 
 \end{proof}

Next, let us turn our attention to the most difficult cases $c_1(\E)=1,2,3$, 
where we have no clear indications on the degree of $C$. We are 
able to dispose of these cases only through a boring list of computations.

\begin{case} There are no indecomposable rank 2 ACM bundles $\E$ on a general sextic, with $c_1(\E)=1$.
\end{case}
\begin{proof} We know from section 3 that we have $\deg(C)=14-h^0(\Ii_C(2))$. For these subcanonical curves in $\Pj^4$, we have yet computed the dimension of the Hilbert scheme in our paper \cite{CM} (see theorem 1.3 and its proof). $h^0(\Ii_C(6))$ is easy to compute since $C$ is ACM, so $\Oc_C(6)$ is non--special. Recall that by \cite{CM} proposition 4.11, we must have $\deg(C)\geq 11$.
The computation are resumed in the following table:
$$\begin{matrix}
\deg(C) & h^0(N(C)) & h^0(\Ii(6)) & h^0(N(C)) + h^0(\Ii(6))-1\\
14	&	56	&	140 	& 195 	\\ 
13	&	53	&	145 	& 197 	\\ 
12	&	50	&	150 	& 199 	\\ 
11	&	47	&	155 	& 201 	\end{matrix}
$$
In any event, the last column is smaller than $209$, so that the map $I\to\Pj$ cannot be dominant
\end{proof}

\begin{case} There are no indecomposable rank 2 ACM bundles $\E$ on a general sextic, with $c_1(\E)=2$.
\end{case}
\begin{proof} We know from section 3 that we have $\deg(C)=20-2h^0(\Ii_C(2))$. 
We have $h^3(\E)=h^0(\E(-1))=0$, hence $\E(3)$ is regular (in the sense of Castelnuovo--Mumford). It follows from sequence (0) that the ideal of $C$ is generated in degree $\leq 5$ (we cannot say that $C$ lies in a smooth quintic, because $C$ itself may be singular).
The resolution of $\Ii$ is of the following type:
$$
\begin{aligned}
0  & \to \Oc(-8)\to \Oc(-6)^c \oplus \Oc(-5)^a \oplus \Oc(-4)^b \oplus \Oc(-3)^x  \to \\
& \to \Oc(-5)^x \oplus \Oc(-4)^b \oplus \Oc(-3)^a \oplus \Oc(-2)^c \to \Ii\to 0 
\end{aligned}
$$

\begin{lemma} Assume that $C$ is contained in $c\geq 3$ independent quadrics. 
Then in fact $c=3$, $\deg(C)=14$ and the resolution of the ideal $\Ii$ of $C$ in $\Pj^4$ is:
$$
\begin{aligned}
0\to \Oc(-8) & \to \Oc(-6)^3 \oplus \Oc(-3)^2 \to \\ 
 & \to \Oc(-5)^2\oplus \Oc(-2)^3 \to \Ii\to 0. 
\end{aligned}
$$ 
\end{lemma}
\begin{proof} First observe that the quadrics through $C$ cannot have a common 
hyperplane. Indeed otherwise, since $C$ is non--degenerate, then one of its components lies outside this hyperplane, and it is a line (remind that $C$ is locally complete intersection, so it has no zero--dimensional components). We get a contradiction, since we yet observed that a general sextic hypersurface contains no lines.\par\noindent
It follows that two general quadrics through $C$ meet along a surface. So, since $C$ is cut scheme-theoretically by quintics, we may link $C$ 
to a curve $C'$ in a complete intersection of type 
$(2,2,5)$. We may compute a minimal resolution of the ideal sheaf of $C'$, using the mapping cone of the diagram:
$$
\begin{matrix}
0  & \to & \Oc(-9)&\to& \Oc(-4)\oplus \Oc(-7)^2 &  \to & \Oc(-5)\oplus \Oc(-2)^2 & \to & \Ii_{2,2,5} \cr
  &  & \downarrow & & \downarrow & & \downarrow &  &   \cr
0  & \to & \Oc(-8)&\to& \Oc(-6)^c \oplus \Oc(-5)^a \oplus &\to & \Oc(-5)^x \oplus \Oc(-4)^b \oplus & \to & \Ii \cr
  &  & & & \oplus \Oc(-4)^b \oplus \Oc(-3)^x &   &  \oplus \Oc(-3)^a \oplus \Oc(-2)^c & &   
\end{matrix}
$$
It turns out that $C'$ is degenerate and its ideal has a minimal syzygy of degree $7$ (it comes from the third quadric containing $C$).\par\noindent
Link now again $C'$ with a complete intersection of type $(1,2,5)$, and use again the mapping cone to compute a resolution of the ideal sheaf of the new curve $C''$. $C''$ is contained in a plane and it has even degree $\leq 4$. Now, we know everything about the resolution of the ideal of $C''$ and we may use this information to compute back the resolution of $C$. Imposing auto--duality, one gets immediately that the degree of $C''$ is $4$ and $x=2$, $a=b=0$. The claim follows.
\end{proof}

\begin{remark} \rm Of course, (smooth) curves with a resolution as the one quoted in the 
previous statement can be easily produced, with a biliaison, starting with a plane quartic curve. \end{remark}

We concluded that the minimal degree for $C$ is $14$.\par\noindent
When the degree is minimal, we know enough of the minimal resolution 
to compute $H^0(N(C))$, using formula (3). Indeed one gets $h^0(N(C))=62$. On the other hand $h^0(\Ii(6))=147$ and $\dim(I)\leq 208$: $I$ cannot dominate $\Pj$ (but it was close!).\par\noindent
Assume $\deg(C)=16$, so that $C$ lies in just two independent quadrics. Then we may have at most one cubic syzygy and the resolution of $\Ii$ now reads:
$$
\begin{aligned}
0  & \to \Oc(-8)\to \Oc(-6)^2 \oplus \Oc(-5)^a \oplus \Oc(-4)^b \oplus \Oc(-3)^x  \to \\
& \to \Oc(-5)^x \oplus \Oc(-4)^b \oplus \Oc(-3)^a \oplus \Oc(-2)^2 \to \Ii\to 0 
\end{aligned}
$$ 
with $x=0,1$. In any event, taking the first Chern class in the sequence, one gets $a=x$ and the computation of $h^0(N(C))$ leads to the same result $66$, whatever $a,b$ are. Since $h^0(\Ii(6))=138$, one concludes as above.\par\noindent
If $\deg(C)=18$, then $C$ lies in one quadric and there are no cubic syzygies. The resolution is:
$$
\begin{aligned}
0 & \to \Oc(-8)\to \Oc(-6)\oplus\Oc(-5)^2 \oplus \Oc(-4)^b \to \\
& \to \Oc(-4)^b \oplus \Oc(-3)^2 \oplus \Oc(-2)\to \Ii\to 0
\end{aligned}
$$ 
($h^0(\Ii_C(3))$ can be computed directly). Then $h^0(N(C))=70$ while $h^0(\Ii_C(6))=129$ 
and one concludes.\par\noindent
Finally, when $\deg(C)=20$, we have no quadrics through $C$ and the resolution reads
$$
0\to \Oc(-8) \to \Oc(-5)^4 \oplus \Oc(-4)^b \to \Oc(-4)^b \oplus \Oc(-3)^4\to \Ii\to 0 
$$ 
from which $h^0(N(C))=74$; since $h^0(\Ii(6))=120$, 
one concludes the proof of this case.
\end{proof}

We arrive now at the (by far) most difficult case $c_1=3$; we deal with $4$--subcanonical curves $C$.

\begin{case} There are no indecomposable rank 2 ACM bundles $\E$ on a general sextic, with $c_1(\E)=3$.
\end{case}

\begin{proof} We know from the previous section that we have $\deg(C)=30-h^0(\Ii_C(3))$. Furthermore $h^3(\E(-1))=h^0(\E(-1))=0$, hence $\E(2)$ is regular (in the sense of Castelnuovo--Mumford). It follows from sequence (0) that the ideal of $C$ is generated by quintics. Moreover, since $\E$ is normalized then $C$ lies in no quadrics.\par\noindent
Call $x$ the number of independent cubics containing $C$ ($x=30-\deg(C)$). Since minimal generators of degree $6$ give, by auto--duality, syzygies of degree $3$, we have none of them. The resolution is of the form:
$$
\begin{aligned}
0 & \to \Oc(-9)\to \Oc(-6)^x\oplus\Oc(-5)^a\oplus\Oc(-4)^b \to \\
& \to \Oc(-5)^b\oplus\Oc(-4)^a\oplus
\Oc(-3)^x \to \Ii\to 0 
\end{aligned}
$$ 
By Riemann-Roch one gets $a+3(27-\deg(C))=b$. Performing the (boring but easy) computations, it turns out that $h^0(N(C))=69+\deg(C)$ while $h^0(\Ii(6))=210-4\deg(C)$, so that  $\dim I= 209+69-3\deg(C)$. To end the proof, it is enough to show that the degree of $C$ cannot be less than $24$.\par\noindent

\begin{claim} Assume $\deg(C)<24$. Then the cubics through $C$ have a common component.
\end{claim}
\begin{proof}
Indeed assume that the cubics through $C$ have no common components. Taking two general such cubics $F_1, F_2$, we can link our curve in a complete intersection of 
type $(3,3,5)$ (since it is cut by quintics). On the other hand, computing the 
mapping cone, it turns out that the linked curve $C'$ has just one minimal 
generator of degree $2$ and it is cut by quintics. A resolution of the ideal sheaf of $C'$ reads:
$$
\begin{aligned}
0  \to &\Oc(-8)^{x-2}\oplus\Oc(-7)^a\oplus \Oc(-6)^{b-1}\to \\ 
\to &\Oc(-7)^b\oplus\Oc(-6)^a\oplus\Oc(-5)^x \to \Oc(-5)\oplus\Oc(-3)^2\oplus\Oc(-2)\to \Ii(C')\end{aligned}
$$ 
The quadric through $C'$ cannot have a common component with the general cubic in the system generated by $F_1$ and $F_2$; so we may link again $C'$ with a complete intersection of type $(2,3,5)$ and get a new curve $C''$ of degree $\leq 8$.
\par\noindent
We claim that $C''$ is cut by quartics. Indeed observe that the ideal sheaf of $C'$ has no syzygies of degree less than $5$, hence all the trivial syzygies between the quadric and one cubic generator of $C'$ are minimal syzygies of the resolution. It follows that, in the mapping cone of the linkage $C'\sim C''$, the trivial syzygy of the cubic and the quadric splits, so that the ideal sheaf of $C''$ has resolution of type:
 $$
\begin{aligned}
0 \to &\Oc(-7)\to \Oc(-5)^{x-1}\oplus\Oc(-4)^a\oplus \Oc(-3)^{b+\epsilon-1}\to \\
\\ &\to\Oc(-4)^{b-1}\oplus\Oc(-3)^{a+\epsilon}\oplus\Oc(-2)^{x-1} \to \Ii(C'')
\end{aligned}
$$ 
where $\epsilon=0,1$. We have $x-1\geq 6$, so the quadrics through $C''$ cannot have a common hyperplane, because they are too many. Thus we may link $C''$ with a complete intersection of type $(2,2,4)$. The curve $C'''$ that we get lies in a hyperplane, has two minimal generators of degree $2$ and it is cut by quartics. Since $\deg(C''')\geq 8$ and $C'''$ is not complete intersection of type $(1,2,4)$, we obtain a contradiction.
\end{proof}

It follows that, in the previous situation,  the cubics through $C$ have 
a fixed component $H$. $H$ is a hyperplane, not a quadric, for when 
$\deg(C)<24$ we have at least $7$ independent cubics through $C$.\par\noindent
Since $C$ is not contained in quadrics, then it has a degenerate part $C_1\subset H$, 
of degree $d_1$; call $C_2$ the residue and $d_2=\deg(C)-d_1$.\par\noindent

\begin{claim} The set--theoretic intersection of the $7$ independent quadrics which contain $C_2$ is a surface $S$.\end{claim}
\begin{proof}
Clearly the quadrics cannot have a common component. If one irreducible component of the intersection is a curve $\Gamma$, then $\Gamma$ has degree $\leq 8$. A classification of curves of low degree on a general sextic threefold was studied by Wu (\cite{Wu}); it turns out that $\Gamma$ must be reduced of degree $6$ and it must be the intersection of $X$ and a plane $\pi$. This implies that all the quadrics which contain $\Gamma$ also contain $\pi$.
\end{proof}

Let us now examine the surface $S$.\par\noindent
The degree of $S$ is not bigger than $4$; it cannot be $4$, for in this case 
$S$ is contained only in two independent quadrics; for the same reason $\deg(S)$ cannot be equal to $3$: in this case $S$ is linked to a plane and one can compute the number of quadrics through it. Similarly $S$ cannot be the union of two distinct planes or an irreducible quadric.\par\noindent
   
Hence either $S$ is a double plane or it is a plane. In any event $(C_2)_{red}$ is a plane curve. By \cite{Wu} again, it must be irreducible of degree $6$ (no curves of lower degree are allowed on a general sextic threefold). Notice that $C_2$ is non--reduced, for $C$ does not lie on quadrics; hence $\deg(C_2)\geq 12$.\par\noindent

Now we arrive at the end of the argument. A general quadric of the system is smooth at a general point of $S$, for we have at most $3$ independent quadrics with vertex at a fixed plane. Taking a general hyperplane section, one sees that $S$ cannot be a double plane, for we do not have a double line on an irreducible quadric surface of $\Pj^3$ which lies in $7$ independent quadrics. It follows that $S$ is 
a plane. \par\noindent
The scheme theoretical intersection of the $7$ quadrics cannot have an embedded curve of degree bigger than $1$. Indeed if $Z$ is a curve in $S$ where all the quadrics are tangent, we may take generators of the linear system such that all but one are singular along $Z$; hence $Z$ is linear. It follows that, in particular,
no components of $C_2$ are embedded parts of this intersection. 
Finally observe that $C_2$ cannot lie in the plane $S$, since $C$ is not contained in quadrics.\par\noindent
This shows that $C_2$ cannot exist: the proof is extablished.
\end{proof}

\begin{remark} \rm Buchweitz, Greuel and Schreyer proved in \cite{B-G-S} that every smooth threefold of degree $>2$ has some indecomposable vector bundle without intermediate cohomology (in fact, it contains a non--discrete family of such bundles). Our main theorem proves that, on a general sextic threefold, such bundles must have rank at least 3 (compare with \cite{Mad4}). 
\end{remark}

\bigskip

\noindent Luca Chiantini \par\noindent
Dipartimento di Scienze Matematiche e Informatiche, Via del Capitano, 15, 53100 SIENA, Italia \par \noindent \it chiantini@unisi.it \rm

\

\noindent Carlo Madonna \par\noindent
Dipartimento di Matematica Univ. Ferrara,
Via Machiavelli, 35 - 44100 Ferrara, Italia \par\noindent
\it madonna@mat.uniroma3.it
\end{document}